\documentclass[a4paper]{article}
\usepackage{amsmath,amssymb,mathrsfs,graphicx,psfrag}
\DeclareMathOperator{\re}{Re}
\DeclareMathOperator{\im}{Im}
\DeclareMathOperator{\fp}{f.p.}
\newtheorem{thm}{Theorem}
\allowdisplaybreaks[1]
\begin{document}
\title{Numerical method of computing Hadamard finite-part integrals with a non-integral power singularity
at the endpoint over a half infinite interval}
\author{Hidenori Ogata\footnote{%
Department of Computer and Network Engineering, 
Graduate School of Informatics and Engineering, 
The University of Electro-Communications, 
1-5-1 Chofugaoka, Chofu, Tokyo 182-8585, Japan, 
(e-mail) {\tt ogata@im.uec.ac.jp}
}}
\maketitle
\begin{abstract}
 In this paper, we propose a numerical method of computing an\\ Hadamard finite-part integral, 
 a finite value assigned to a divergent integral, with a non-integral power singularity at the endpoint 
 on a half infinite interval. 
 In the proposed method, we express a desired finite part integral using a complex integral, 
 and we obtain the finite part integral by evaluating the complex integral by the DE formula. 
 Theoretical error estimate and some numerical examples show the exponential convergence of the proposed method. 
\end{abstract}
\section{Introduction}
\label{sec:introduction}
The integral
\begin{equation*}
 \int_0^{\infty}x^{\alpha-2}f(x)\mathrm{d}x \quad ( \: 0 < \alpha < 1 \: ),
\end{equation*}
where $f(x)$ is an analytic function on $[0,+\infty)$ such that $f(0)\neq 0$ and 
$f(x)=\mathrm{O}(x^{1-\alpha-\varepsilon})$ as $x\rightarrow+\infty$ $( \: \varepsilon>0 \: )$, 
is divergent. 
However, we can assign a finite value to this divergent integral. 
In fact, for $\epsilon>0$, we have by integration by part
\begin{align*}
 \int_{\epsilon}^{\infty}x^{\alpha-2}\mathrm{d}x = \: & 
 -\frac{1}{1-\alpha}\int_{\epsilon}^{\infty}(x^{\alpha-1})^{\prime}f(x)\mathrm{d}x
 \\ 
 = \: & 
 -\frac{1}{1-\alpha}
 \left\{
 \bigg[x^{\alpha-1}f(x)\bigg]_{\epsilon}^{\infty} - \int_{\epsilon}^{\infty}x^{\alpha-1}f^{\prime}(x)\mathrm{d}x
 \right\}
 \\ 
 = \: & 
 \frac{\epsilon^{\alpha-1}f(\epsilon)}{1-\alpha}
 + 
 \frac{1}{1-\alpha}\int_{\epsilon}^{\infty}x^{\alpha-1}f^{\prime}(x)\mathrm{d}x
 \\ 
 = \: & 
 \frac{\epsilon^{\alpha-1}f(0)}{1-\alpha} + \mathrm{O}(1) 
 \quad \mbox{as} \quad \epsilon\downarrow 0, 
\end{align*}
and the limit
\begin{equation*}
 \lim_{\epsilon\downarrow 0}
  \left\{ 
   \int_{\epsilon}^{\infty}x^{\alpha-2}f(x)\mathrm{d}x
   - 
   \frac{\epsilon^{\alpha-1} f(0)}{1-\alpha}
  \right\}
\end{equation*}
exists and is finite. 
We call this limit an Hadamard finite-part (f.p.) integral and denote it by
\begin{equation*}
 \fp\int_0^{\infty}x^{\alpha-2}f(x)\mathrm{d}x.
\end{equation*}
More generally, we can define the f.p. integral
\begin{equation}
 \label{eq:fp-integral0}
  \fp\int_0^{\infty}x^{\alpha-1-n}f(x)\mathrm{d}x
\end{equation}
for $n=1,2,\ldots$, $0<\alpha<1$ and an analytic function $f(x)$ on $[0,+\infty)$ 
such that $f(0)\neq 0$ and $f(x) = \mathrm{O}(x^{n-\alpha-\varepsilon})$ as $x\rightarrow+\infty$ 
$( \: \varepsilon > 0 \: )$ \cite{EstradaKanwal1989}. 

In this paper, we propose a numerical method of computing a f.p. integral (\ref{eq:fp-integral0}). 
In the proposed method, we express the desired f.p. integral using a complex integral, 
and obtain the f.p. integral by evaluating the complex integral by the DE formula 
\cite{TakahasiMori1974}. 
Theoretical error estimate and some numerical examples show that the proposed approximation formula converges 
exponentially as the number of sampling points increases. 

Previous studies related to this paper are as follows. 
The author and Hirayama proposed a numerical method of computing ordinary integrals based on 
hyperfunction theory, a theory of generalized functions based on complex function theory 
\cite{OgataHirayama2018}. 
In their method, we obtain the desired integral by evaluating a complex integral using a conventional numerical 
integration formula as in the method for computing f.p. integrals proposed in this paper. 
The author proposed numerical methods of computing a f.p. integral with a singularity at an endpoint 
on a finite interval \cite{Ogata2019c,Ogata2019b} and a f.p. integral with an integral power singularity 
at the endpoint on a half infinite interval \cite{Ogata2019d}.  
Also in these methods, we obtain a desired integral using a complex integral, and obtain the integral 
by evaluating the complex integral by a conventional numerical integration formula. 
For the computation of a Cauchy principal value integral or 
a f.p. integral on a finite interval with a singularity in the interior 
of the integral interval
\begin{equation}
 \label{eq:fp-integral02}
  \fp\int_a^b \frac{f(x)}{(x-\lambda)^n}\mathrm{d}x 
  \quad ( \: -\infty < a < \lambda < b < +\infty, \: n = 1, 2, \ldots \: ), 
\end{equation}
many numerical methods were proposed. 
Elliot and Paget proposed Gauss-type numerical integration formulas for (\ref{eq:fp-integral02}) 
\cite{ElliotPaget1979,Paget1981}. 
Bialecki proposed Sinc numerical integration formula of computing (\ref{eq:fp-integral02}), 
where the trapezoidal formula together with a variable transform technique are used 
as in the DE formula \cite{TakahasiMori1974}. 
The author et al. improved these methods and proposed a DE-type numerical integration formula of 
computing (\ref{eq:fp-integral02}) \cite{OgataSugiharaMori2000}. 

The remainder of this paper is structured as follows. 
In Section \ref{sec:fp-integral}, we define the f.p. integral (\ref{eq:fp-integral0}) 
and propose a numerical method of computing it. 
In addition, we show a theoretical error estimate which shows the exponential convergence of the proposed method. 
In Section \ref{sec:example}, we show some numerical examples which show the effectiveness of the proposed method. 
In Section \ref{sec:summary}, we give a summary of this paper.
\section{Hadamard finite-part integral and a numerical method}
\label{sec:fp-integral}
Let $n=1,2,\ldots$, $0<\alpha<1$, and $f(x)$ be an analytic function on $[0,+\infty)$ such that 
$f(0)\neq 0$ and $f(x)=\mathrm{O}(x^{n-\alpha-\epsilon})$ as $x\rightarrow+\infty$ $( \: 0 < \varepsilon < 1 \: )$. 
The Hadamard finite-part integral (\ref{eq:fp-integral0}) is defined by
\begin{align}
 \nonumber
 I^{(n,\alpha)}[f] =  \: & 
 \fp\int_0^{\infty}x^{\alpha-1-n}f(x)\mathrm{d}x 
 \\
 \label{eq:fp-integral}
 = \: & 
 \lim_{\epsilon\downarrow 0}
 \left\{
 \int_{\epsilon}^{\infty}x^{\alpha-1-n}f(x)\mathrm{d}x 
 - 
 \sum_{k=0}^{n-1}\frac{\epsilon^{\alpha-n+k}}{k!(n-\alpha-k)}f^{(k)}(0)
 \right\}. 
\end{align}
We can show that (\ref{eq:fp-integral}) is well-defined as follows. 
In fact, by integration by part, we have for $\epsilon>0$
\begin{align*}
 & 
 \int_{\epsilon}^{\infty}x^{\alpha-1-n}f(x)\mathrm{d}x
 \\ 
 = \: & 
 -\frac{1}{n-\alpha}\int_{\epsilon}^{\infty}(x^{\alpha-n})^{\prime}f(x)\mathrm{d}x
 \\ 
 = \: & 
 - \frac{1}{n-\alpha}
 \left\{
 \bigg[x^{\alpha-n}f(x)\bigg]_{\epsilon}^{\infty}
 - 
 \int_{\epsilon}^{\infty}x^{\alpha-n}f^{\prime}(x)\mathrm{d}x
 \right\}
 \\ 
 = \: & 
 \frac{\epsilon^{\alpha-n}f(\epsilon)}{n-\alpha} 
 + 
 \frac{1}{n-\alpha}\int_{\epsilon}^{\infty}x^{\alpha-n}f^{\prime}(x)\mathrm{d}x
 \\ 
 = \: & 
 \frac{\epsilon^{\alpha-n}f(\epsilon)}{n-\alpha}
 + 
 \frac{\epsilon^{\alpha-n+1}f^{\prime}(\epsilon)}{(n-\alpha)(n-\alpha-1)}
 + 
 \frac{1}{(n-\alpha)(n-\alpha-1)}\int_{\epsilon}^{\infty}
 x^{\alpha-n+1}f^{\prime\prime}(x)\mathrm{d}x
 \\
 = \: & \cdots 
 \\ 
 = \: & 
 \sum_{k=0}^{n-1}\frac{\epsilon^{\alpha-n+k}f^{(k)}(\epsilon)}{(n-\alpha)(n-\alpha-1)\cdots(n-\alpha-k)}
 \\
 & 
 + 
 \frac{1}{(n-\alpha)(n-\alpha-1)\cdots(1-\alpha)}\int_{\epsilon}^{\infty}x^{\alpha-1}f^{(n)}(x)\mathrm{d}x
  \\ 
 = \: & 
 \sum_{k=0}^{n-1}\frac{\epsilon^{\alpha-n+k}}{(n-\alpha)(n-\alpha-1)\cdots(n-\alpha-k)}
 \sum_{l=0}^{\infty}\frac{\epsilon^l}{l!}f^{(l+k)}(0)
 \\
 & 
 + 
 \frac{1}{(n-\alpha)(n-\alpha-1)\cdots(1-\alpha)}\int_{\epsilon}^{\infty}x^{\alpha-1}f^{(n)}(x)\mathrm{d}x
 \\
 = \: & 
 \sum_{k=0}^{n-1}\sum_{l=0}^{n-k-1}
 \frac{\epsilon^{\alpha-n+k+l}f^{(l+k)}(0)}{l!(n-\alpha)(n-\alpha-1)\cdots(n-\alpha-k)} + \mathrm{O}(1)
 \\
 = \: & 
 \sum_{m=0}^{n-1}\left\{
 \sum_{k+l=m}
 \frac{1}{l!(n-\alpha)(n-\alpha-1)\cdots(n-\alpha-k)} 
 \right\}
 \epsilon^{\alpha-n+m}f^{(m)}(0) 
 + \mathrm{O}(1)
 \\
 = \: & 
 \sum_{m=0}^{n-1}
 \underbrace{%
 \left\{
 \sum_{k=0}^{m}\frac{1}{(m-k)!(n-\alpha)(n-\alpha-1)\cdots(n-\alpha-k)}
 \right\}
 }_{(\ast)}
 \epsilon^{\alpha-n+m}f^{(m)}(0)
 \\
 & 
 + \mathrm{O}(1) \quad \mbox{as} \ \epsilon\downarrow 0, 
\end{align*}
and 
\begin{align*}
 (\ast) = \: & 
 \frac{1}{m!(n-\alpha)} + \frac{1}{(m-1)!(n-\alpha)(n-\alpha-1)} + \cdots
 \\ 
 & 
 + 
 \frac{1}{1!(n-\alpha)\cdots(n-\alpha-m+2)(n-\alpha-m+1)}
 \\
 & 
 +
 \frac{1}{(n-\alpha)\cdots(n-\alpha-m+2)(n-\alpha-m+1)(n-\alpha+m)}
 \\
 = \: & 
 \frac{1}{m!(n-\alpha)} + \frac{1}{(m-1)!(n-\alpha)(n-\alpha-1)} + \cdots
 \\
 & 
 + 
 \frac{1}{2!(n-\alpha)\cdots(n-\alpha-m+3)(n-\alpha-m+2)}
 \\
 & 
 + 
 \frac{1}{1!(n-\alpha)\cdots(n-\alpha-m+3)(n-\alpha-m+2)(n-\alpha-m)}
 \\
 = \: & \cdots 
 = 
 \frac{1}{m!(n-\alpha-m)}.
\end{align*}
Then, we have
\begin{equation*}
 \int_{\epsilon}^{\infty}x^{\alpha-1-n}f(x)\mathrm{d}x = 
  \sum_{m=0}^{n-1}\frac{\epsilon^{\alpha-n+m}}{m!(n-\alpha-m)}f^{(m)}(0)
  +
  \mathrm{O}(1) \quad \mbox{as} \ \epsilon\downarrow 0,
\end{equation*}
and the limit of (\ref{eq:fp-integral}) exists and is finite. 

\medskip

\begin{figure}[htbp]
 \begin{center}
  \psfrag{0}{$\mathrm{O}$}
  \psfrag{e}{$\epsilon$}
  \psfrag{p}{$\Gamma_{\epsilon}^{(+)}$}
  \psfrag{m}{$\Gamma_{\epsilon}^{(-)}$}
  \psfrag{a}{$C_{\epsilon}$}
  \includegraphics[width=0.7\textwidth]{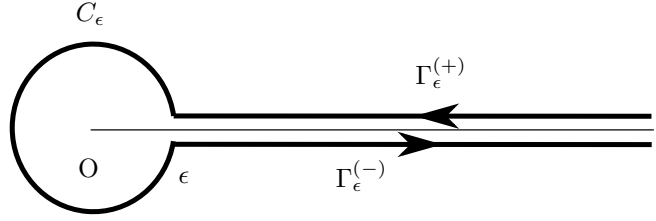}
 \end{center}
 \caption{The complex integral paths $\Gamma_{\epsilon}^{(\pm)}$ and $C_{\epsilon}$.}
 \label{fig:proof-integral-path}
\end{figure}

As shown in the following theorem, a f.p. integral (\ref{eq:fp-integral}) is expressed using a complex integral, 
which is the bases of our numerical method. 
\begin{thm}
 \label{thm:complex-integral}
 We suppose that $f(z)$ is analytic in a domain $D$ containing the positive real axis in its interior. 
 Then, we have
 \begin{equation}
  \label{eq:complex-integral}
   I^{(n,\alpha)}[f]
   = 
   \frac{(-1)^{n+1}}{2\mathrm{i}\sin\pi\alpha}
   \oint_C (-z)^{\alpha-1-n}f(z)\mathrm{d}z, 
 \end{equation}
 where $C$ is a complex integral path such that it encircles the positive real axis in the positive sense 
 and it is contained in $D$, and $z^{\alpha-1-n}$ is the principal value, that is, the branch such that 
 it takes a real value on the positive real axis%
 \footnote{%
 The complex integral on the right-hand side of (\ref{eq:complex-integral}) coincides with %
 the integral of $x^{\alpha-1-n}f(x)$ as a hyperfunction \cite{Graf2010}. %
 }. 
\end{thm}
\paragraph{Proof of Theorem \ref{thm:complex-integral}}
By Cauchy's integral theorem, we have
\begin{equation*}
 \oint_C (-z)^{\alpha-1-n}f(z)\mathrm{d}z = 
  \left( \int_{C_{\epsilon}} + \int_{\Gamma_{\epsilon}^{(+)}} + \int_{\Gamma_{\epsilon}^{(-)}}\right)
  (-z)^{\alpha-1-n}f(z)\mathrm{d}z,
\end{equation*}
where $\epsilon>0$, and 
$\Gamma_{\epsilon}^{(\pm)}$ and $C_{\epsilon}$ are the complex integral paths respectively given by
\begin{align*}
 \Gamma_{\epsilon}^{(+)} = \: & 
 \{ \: x + \mathrm{i}0 \: | \: +\infty > x \geqq \epsilon \: \}, 
 \\ 
 \Gamma_{\epsilon}^{(-)} = \: & 
 \{ \: x - \mathrm{i}0 \: | \: \epsilon \leqq x < +\infty \: \}, 
 \\
 C_{\epsilon} = \: & 
 \{ \: \epsilon\mathrm{e}^{\mathrm{i}\theta} \: | \: 0 \leqq \theta \leqq 2\pi \: \}
\end{align*}
(see Figure \ref{fig:proof-integral-path}).
Regarding the integrals on $\Gamma_{\epsilon}^{(\pm)}$, we have
\begin{align*}
 & 
 \left( \int_{\Gamma_{\epsilon}^{(+)}} + \int_{\Gamma_{\epsilon}^{(-)}}\right)
 (-z)^{\alpha-1-n}f(z)\mathrm{d}z
 \\
 = \: & 
 - \int_{\epsilon}^{\infty}(-(x+\mathrm{i}0))^{\alpha-1-n}f(x)\mathrm{d}x
 + \int_{\epsilon}^{\infty}(-(x-\mathrm{i}0))^{\alpha-1-n}f(x)\mathrm{d}x
 \\
 = \: & 
 - \int_{\epsilon}^{\infty}(x\mathrm{e}^{-\mathrm{i}\pi})^{\alpha-1-n}f(x)\mathrm{d}x
 + \int_{\epsilon}^{\infty}(x\mathrm{e}^{\mathrm{i}\pi})^{\alpha-1-n}f(x)\mathrm{d}x
 \\ 
 = \: & 
 2\mathrm{i}(-1)^{n+1}\sin\pi\alpha
 \int_{\epsilon}^{\infty}x^{\alpha-1-n}f(x)\mathrm{d}x.
\end{align*}
Regarding the integral on $C_{\epsilon}$, we have
\begin{align*}
 & 
 \int_{C_{\epsilon}}(-z)^{\alpha-1-n}f(z)\mathrm{d}z
 \\
 = \: & 
 \int_0^{2\pi}(\epsilon\mathrm{e}^{\mathrm{i}(\theta-\pi)})^{\alpha-1-n}
 f(\epsilon\mathrm{e}^{\mathrm{i}\theta})\mathrm{i}\epsilon\mathrm{e}^{\mathrm{i}\theta}
 \mathrm{d}\theta
 \\
 = \: & 
 \mathrm{i}(-1)^{n+1}\epsilon^{\alpha-n}\mathrm{e}^{-\mathrm{i}\pi\alpha}
 \int_0^{2\pi}\mathrm{e}^{\mathrm{i}(\alpha-n)\theta}f(\epsilon\mathrm{e}^{\mathrm{i}\theta})
 \mathrm{d}\theta
 \\
 = \: & 
 \mathrm{i}(-1)^{n+1}\epsilon^{\alpha-n}\mathrm{e}^{-\mathrm{i}\pi\alpha}
 \int_0^{2\pi}\mathrm{e}^{\mathrm{i}(\alpha-n)\theta}
 \sum_{k=0}^{\infty}\frac{\epsilon^k}{k!}f^{(k)}(0)\mathrm{e}^{\mathrm{i}k\theta}
 \mathrm{d}\theta
 \\
 = \: & 
 \mathrm{i}(-1)^{n+1}\mathrm{e}^{-\mathrm{i}\pi\alpha}
 \sum_{k=0}^{\infty}\frac{\epsilon^{\alpha-n+k}}{k!}f^{(k)}(0)
 \int_0^{2\pi}\mathrm{e}^{\mathrm{i}(\alpha-n+k)\theta}\mathrm{d}\theta
 \\
 = \: & 
 - 2\mathrm{i}(-1)^{n+1}\sin\pi\alpha
 \sum_{k=0}^{n-1}\frac{\epsilon^{\alpha-n+k}}{k!(n-\alpha-k)}f^{(k)}(0)
 + \mathrm{o}(1) \quad \mbox{as} \ \epsilon\downarrow 0, 
\end{align*}
where we exchanged the order of the integration and the infinite sum on the fourth equality 
since the infinite sum is uniformly convergent on $0\leqq\theta\leqq 2\pi$. 
Summarizing the above calculations, we have
\begin{align*}
 & 
 \oint_C (-z)^{\alpha-1-n}f(z)\mathrm{d}z
 \\
 = \: & 
 2\mathrm{i}(-1)^{n+1}\sin\pi\alpha
 \left\{
 \int_{\epsilon}^{\infty}x^{\alpha-1-n}f(x)\mathrm{d}x 
 - 
 \sum_{k=0}^{n-1}\frac{\epsilon^{\alpha-n+k}}{k!(n-\alpha-k)}f^{(k)}(0)
 \right\}
 \\
 & 
 +
 \mathrm{o}(1) \quad \mbox{as} \ \epsilon\downarrow 0,
\end{align*}
and, taking the limit $\epsilon\downarrow 0$, we obtain (\ref{eq:complex-integral}).
\hfill\rule{1.5ex}{1.5ex}

\medskip

We obtain an approximation formula of computing the f.p. integral (\ref{eq:fp-integral}) by 
evaluating the complex integral in (\ref{eq:complex-integral}) by the DE formula \cite{TakahasiMori1974}
\begin{align}
 \nonumber 
 \int_{-\infty}^{\infty}g(u)\mathrm{d}u = \: & 
 \int_{-\infty}^{\infty}g(\psi_{\rm DE}(v))\psi_{\rm DE}^{\prime}(v)\mathrm{d}v
 \\ 
 \simeq \: & 
 h\sum_{k=-N_-}^{N_+}g(\psi_{\rm DE}(kh))\psi_{\rm DE}^{\prime}(kh), 
\end{align}
where $\psi_{\rm DE}$ is the DE transform
\begin{equation}
 \psi_{\rm DE}(v) = 
  \begin{cases}
   \sinh(\sinh v) & \mbox{if} \ g(u) = \mathrm{O}(|u|^{-1-\alpha}) \quad \mbox{as} \ u\rightarrow\pm\infty \ 
   ( \: \alpha > 0 \: ) 
   \\
   \sinh v & \mbox{if} \ g(u) = \exp(-c|u|) \quad \mbox{as} \ u\rightarrow\pm\infty \ ( \: c > 0 \: ).
  \end{cases}
\end{equation}
We can take the positive integers $N_{\pm}$ small for a given mesh $h$ since the transformed integrand 
$g(\psi_{\rm DE}(v))\psi_{\rm DE}^{\prime}(v)$ decays double exponentially as $u\rightarrow\pm\infty$. 
Taking a parameterization of the integral path 
\begin{equation*}
 C \: : \: z = \varphi(u), \quad - \infty < u < +\infty,
\end{equation*}
and evaluating the complex integral of (\ref{eq:complex-integral}) by the DE formula, 
we obtain the approximation formula of the f.p. integral
\begin{align}
 \nonumber 
 I^{(n,\alpha)}[f] = \: & 
 \frac{(-1)^{n+1}}{2\mathrm{i}\sin\pi\alpha}
 \int_{-\infty}^{\infty}(-\varphi(u))^{\alpha-1-n}f(\varphi(u))\varphi^{\prime}(u)\mathrm{d}u
 \\ 
 \nonumber 
 \simeq \: & I_{h,N_+,N_-}^{(n,\alpha)}[f]
 \\ 
 \nonumber
 = \: & 
 \frac{(-1)^{n+1}h}{2\mathrm{i}\sin\pi\alpha}\sum_{k=-N_-}^{N_+}
 (-\varphi(\psi_{\rm DE}(kh)))^{\alpha-1-n}f(\varphi(\psi_{\rm DE}(kh)))
 \\
 \label{eq:approx-fp-integral}
 & \hspace{25mm}
 \times\varphi^{\prime}(\psi_{\rm DE}(kh))\psi_{\rm DE}^{\prime}(kh).
\end{align}

A theoretical error estimate of the approximation (\ref{eq:approx-fp-integral}) is given 
in the following theorem, where $N_{\pm}$ are taken to be $N_+=N_-=N$ for the simplicity. 
\begin{thm}
 \label{thm:error-estimate}
 We suppose that 
 \begin{itemize}
  \item $\varphi(\psi_{\rm DE}(w))$ is an analytic function in the strip
		\begin{equation*}
		 \mathscr{D}_d = 
		  \{ \: w\in\mathbb{C} \: | \: |\im w| < d \: \} \quad ( \: d > 0 \: ),
		\end{equation*}
		and the domain 
		\begin{equation*}
		 \varphi(\psi_{\rm DE}(\mathscr{D}_d)) = 
		  \{ \: \varphi(\psi_{\rm DE}(w))\: | \: w\in\mathscr{D}_d \: \}
		\end{equation*}
		is contained in $\mathbb{C}\setminus [0,+\infty)$,
  \item $f(\varphi(\psi_{\rm DE}(w)))\varphi^{\prime}(\psi_{\rm DE}(w))\psi_{\rm DE}^{\prime}(w)$ satisfies
		\begin{align*}
		 & 
		 \mathscr{N}^{(n,\alpha)}(f,\varphi,\psi_{\rm D}, \mathscr{D}_d) 
		 \\
		 \equiv \: & 
		 \lim_{\epsilon\downarrow 0}\oint_{\mathscr{D}_d(\epsilon)}
		 \left|
		 (-\varphi(\psi_{\rm DE}(w)))^{\alpha-1-n}f(\varphi(\psi_{\rm DE}(w)))
		 \varphi^{\prime}(\psi_{\rm DE}(w))\psi_{\rm DE}(w)
		 \right|
		 |\mathrm{d}w|
		 \\
		 < \: & \infty,
		\end{align*}
		where
		\begin{equation*}
		 \mathscr{D}_d(\epsilon) = 
		  \{ \: w\in\mathbb{C} \: | \: |\re w| < 1/\epsilon, \quad |\im w| < d(1-\epsilon)\:\},
		\end{equation*}
		and
  \item there exists positive numbers $C, c_1$ and $c_2$ such that
		\begin{multline*}
		 \left|(-\varphi(\psi_{\rm DE}(v)))^{\alpha-1-n}
		  f(\varphi(\psi_{\rm DE}(u)))\varphi^{\prime}(\psi_{\rm DE}(v))
		  \psi_{\rm DE}^{\prime}(v)
		 \right|
		 \\
		 \leqq C\exp(-c_1\exp(c_2|v|)) \quad ( \: \forall v\in\mathbb{R} \: ).
		\end{multline*}
 \end{itemize}
 Then, we have the inequality
 \begin{align}
  \nonumber
  & 
  |I^{(n,\alpha)}[f] - I_{h,N}^{(n,\alpha)}[f]| 
  \\
  \nonumber
  \leqq \: & 
  \frac{1}{2\pi}\mathscr{N}^(n,\alpha)(f,\varphi,\psi_{\rm DE},\mathscr{D}_d)
  \frac{\exp(-2\pi d/h)}{1-\exp(-2\pi d/h)}
  \\
  & 
  \label{eq:error-estimate}
  +
  C^{(n,\alpha)}(f,\varphi,\psi_{\rm DE}, \mathscr{D}_d)\exp(-c_1\exp(c_2 Nh)),
 \end{align}
 where 
 \begin{equation*}
  I_{h,N}^{(n,\alpha)}[f] = 
   I_{h,N,N}^{(n,\alpha)}[f]
 \end{equation*}
 and 
 $C^{(n,\alpha)}(f,\varphi,\psi_{\rm DE}, \mathscr{D}_d)$ 
 is a positive number depending on 
 $n, \alpha, f(z), \varphi, \psi_{\rm DE}$ and $\mathscr{D}_d$ only. 
\end{thm}
This theorem shows that the proposed approximation (\ref{eq:approx-fp-integral}) converges exponentially 
as the mesh $h$ decreases and the number of sampling points $2N^{\prime}+1$ increases. 

\paragraph{Proof of Theorem \ref{thm:error-estimate}}
We have
\begin{align}
 \nonumber
 & 
 |I^{(n,\alpha)}[f]-I_{h,N}^{(n,\alpha)}[f]| 
 \\
 \nonumber
 \leqq \: & 
 |I^{(n,\alpha)}[f]-I_{h}^{(n,\alpha)}[f]|
 \\
 \label{eq:proof-error-estimate}
 & + 
 \left|
 \frac{h}{2\pi\mathrm{i}}\sum_{|k|>N}
 (-\varphi(\psi_{\rm DE}(kh)))^{\alpha-1-n}f(\varphi(\psi_{\rm DE}(kh)))
 \varphi^{\prime}(\psi_{\rm DE}(kh))\psi_{\rm DE}(kh)
 \right|,
\end{align}
where $I_{h}^{(n,\alpha)}[f] = \lim_{N\rightarrow\infty}I_{h,N}^{(n,\alpha)}[f]$. 
Regarding the first term on the right-hand side of (\ref{eq:proof-error-estimate}), 
from (\ref{eq:complex-integral}) and Theorem 3.2.1 of \cite{Stenger1993}, we have
\begin{equation*}
 |\mbox{the first term}| \leqq 
  \frac{1}{2\pi}\mathscr{N}(f,\varphi,\psi_{\rm DE},\mathscr{D}_d)
  \frac{\exp(-2\pi d/h)}{1-\exp(-2\pi d/h)}.
\end{equation*}
Regarding the second term, we have
\begin{align*}
 |\mbox{the second term}| \leqq \: & 
 Ch\sum_{|k|>N}\exp(-c_1\exp(c_2 kh))
 \\ 
 \leqq \: & 
 2C\int_{Nh}^{\infty}\exp(-c_1\exp(c_2 x))\mathrm{d}x
 \\
 \leqq \: & 
 2C\int_{Nh}^{\infty}\exp(c_2 x)\exp(-c_1\exp(c_2 x))\mathrm{d}x
 \\
 = \: & 
 \frac{2C}{c_1 c_2}\exp(-c_1\exp(c_2 Nh)).
\end{align*}
Therefore, we have (\ref{eq:error-estimate}).
\hfill\rule{1.5ex}{1.5ex}

\medskip

We remark here that, if $f(z)$ is real valued on the real axis, we can reduce the number of sampling points 
by half. 
In fact, in this case, we have
\begin{equation*}
 f(\overline{z}) = \overline{f(z)}
\end{equation*}
from the reflection principle. 
Then, taking the integral path $C$ to be symmetric with respect to the real axis, that is, 
\begin{equation*}
 \varphi(-u) = \overline{\varphi(u)},
\end{equation*}
which leads to $\varphi^{\prime}(-u) = - \overline{\varphi^{\prime}(u)}$, we have
\begin{align}
 \nonumber
 I^{(n,\alpha)}[f] \simeq \: & I_{h,N}^{(n,\alpha)}[f]
 \\ 
 \nonumber
 = \: & 
 \frac{h}{\sin\pi\alpha}\im\bigg\{
 \frac{1}{2}(-\varphi(\psi_{\rm DE}(0)))^{\alpha-1-n}f(\varphi(\psi_{\rm DE}(0)))
 \varphi^{\prime}(\psi_{\rm DE}(0))\psi_{\rm DE}^{\prime}(0)
 \\
 \label{eq:approx-fp-integral2}
 & + 
 \sum_{k=1}^{N}
 (-\varphi(\psi_{\rm DE}(kh)))^{\alpha-1-n}f(\varphi(\psi_{\rm DE}(kh)))
 \varphi^{\prime}(\psi_{\rm DE}(kh))\psi_{\rm DE}^{\prime}(kh)
 \bigg\}
\end{align}
\section{Numerical examples}
\label{sec:example}
\begin{figure}[htbp]
 \begin{center}
  \psfrag{x}{$\re z$}
  \psfrag{y}{$\im z$}
  \includegraphics[width=0.5\textwidth]{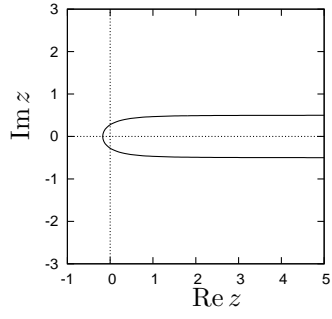}
 \end{center}
 \caption{The complex integral path $C$.}
 \label{fig:example-integral-path}
\end{figure}

\begin{figure}[htbp]
 \begin{center}
  \begin{tabular}{cc}
   \psfrag{N}{$N$}
   \psfrag{e}{\hspace{-10mm}relative error}
   \includegraphics[width=0.45\textwidth]{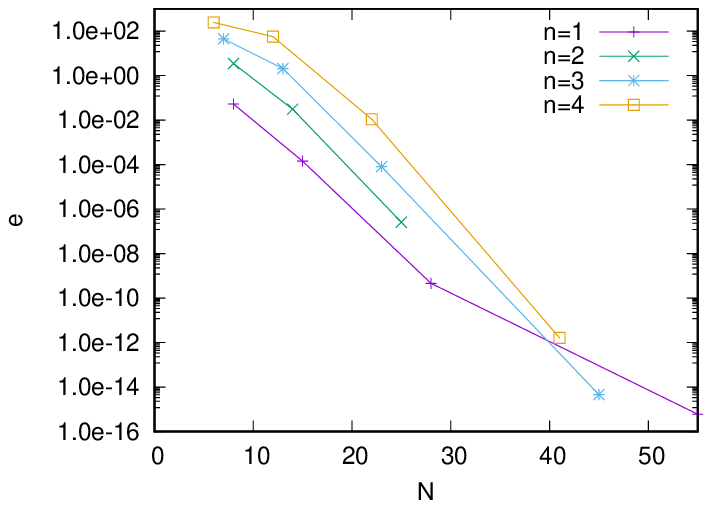}
   & 
       \psfrag{N}{$N$}
       \psfrag{e}{\hspace{-10mm}relative error}
       \includegraphics[width=0.45\textwidth]{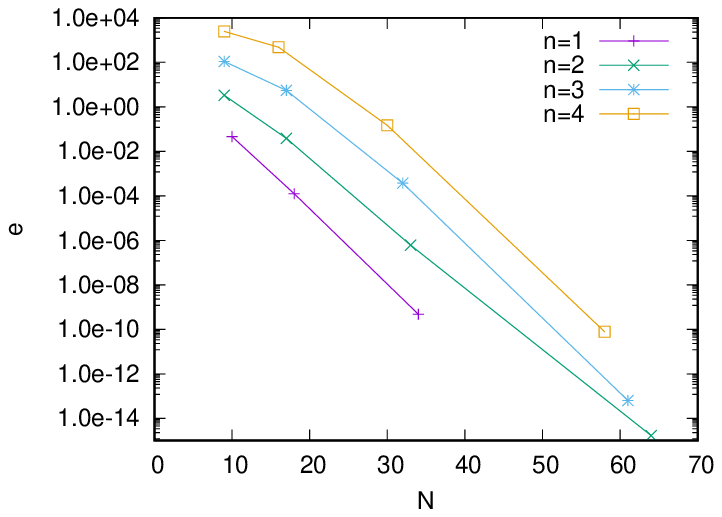}
       \\
   integral (i) & integral (ii)
  \end{tabular}
 \end{center}
 \caption{The errors of the proposed approximation formula (\ref{eq:approx-fp-integral2})
 applied to the f.p. integrals (\ref{eq:example}) 
 with $\alpha = 0.5$.}
 \label{fig:example}
\end{figure}

We computed the f.p. integrals
\begin{equation}
 \label{eq:example}
  \begin{aligned}
   \mathrm{(i)} \quad & 
   \fp\int_0^{\infty}\frac{x^{\alpha-1-n}}{1+x^2}\mathrm{d}x = 
   \begin{cases}
    (-1)^m(\pi/2)/\sin(\pi\alpha/2) & n = 2m \ \mbox{(even)} \\ 
    (-1)^{m+1}(\pi/2)/\cos(\pi\alpha/2) & n = 2m+1 \ \mbox{(odd)},
   \end{cases} 
   \\ 
   \mathrm{(ii)} \quad & 
   \fp\int_0^{\infty}x^{\alpha-1-n}\mathrm{e}^{-x}\mathrm{d}x = 
   \frac{(-1)^n\Gamma(\alpha)}{(1-\alpha)(2-\alpha)\cdots(n-\alpha)}
  \end{aligned}
\end{equation}
with $\alpha = 0.5$ by the proposed approximation formula (\ref{eq:approx-fp-integral2}). 
We performed all the computations using programs coded in C++ with double precision working. 
We took the complex integral path $C$ as
\begin{equation*}
 C \: : \: 
  w = \varphi(u) = 
  \frac{u+0.5\mathrm{i}}{\mathrm{i}\pi}
  \log\left(\frac{1+\mathrm{i}(u+0.5\mathrm{i})}{1-\mathrm{i}(u+0.5\mathrm{i})}\right), 
  \quad 
  + \infty > u > - \infty
\end{equation*}
(see Figure \ref{fig:example-integral-path}).
We decided the number of sampling points $N$ for given mesh $h=2^{-1}, 2^{-2}, \ldots$ by truncating 
the infinite sum of the right-hand side of (\ref{eq:approx-fp-integral2}) 
at the $k$-th term satisfying 
\begin{equation*}
 \frac{h}{\sin\pi\alpha}|\mbox{the transformed integrand of the $k$-th term}| < 
  10^{-15}\times \left|I_{h,N}^{(n,\alpha)}[f]\right|. 
\end{equation*}
Figure \ref{fig:example} shows the relative errors of the proposed method (\ref{eq:approx-fp-integral2})
applied to the f.p. integrals 
(\ref{eq:example}). 
From these figures, the proposed formula converges exponentially as the number of sampling points $N$ increases.
\section{Summary}
\label{sec:summary}
In this paper, we proposed a numerical method of computing a f.p. integral with a non-integral power singularity 
at the endpoint over a half infinite interval. 
In the proposed method, we express the desired f.p. integral by a complex integral, and we obtain the f.p. 
integral by evaluating the complex integral by the DE formula. 
Theoretical error estimate and some numerical examples show that the proposed approximation converges 
exponentially as the mesh of the DE formula decreases and the number of sampling points increases 
for an analytic integrand. 

The complex integral which expresses a desired f.p. integral and gives the basis of the proposed method 
coincides with 
the definition of the integral of a hyperfunction, a generalized function given by an analytic function. 
In hyperfunction theory \cite{Graf2010}, 
a hyperfunction is described by an analytic function called a defining function, 
and its integral is defined by an complex integral involving the defining function. 
In addition, in hyperfunction theory, we can deal with an ordinary integral and a f.p. integral in a unified way.  
We found the proposed method, that is, the method of computing a f.p. integral by evaluating a complex integral, 
from this viewpoint in hyperfunction theory. 
Therefore, we expect that hyperfunction theory is applicable to many numerical computations in science and 
engineering. 

\bibliographystyle{plain}
\bibliography{arxiv2019_5}
\end{document}